**If parallel lines could meet: What exactly can a poet say about the Fano plane?**


Dr Katherine Collins

*Faculty of English Language and Literature*
*University of Oxford, Oxford, UK*

*St Cross Building*
*Manor Road*
*Oxford*
*OX1 3UL*

*katherine.collins@ell.ox.ac.uk*

Dr Siaw-Lynn Ng

*Information Security Group*
*Royal Holloway University of London, UK*

*Egham, Surrey, TW20 0EX*

*s.ng@rhul.ac.uk*


Word count: 6474



# If parallel lines could meet: What exactly can a poet say about the Fano plane?


This article describes our invention of a new poetic form based on projective geometry. In doing this we also explore the 'what ifs' in mathematics and poetry which spark the creative processes of poet and mathematician. In other words, throughout our collaboration we often asked one another, is this what it's like for you? Do you think in this way, too? How does your experience of creativity compare to mine? And often, as well, what exactly do you mean when you say…? We spent a fair amount of time and energy, for example, trying to understand one another's interpretation of 'a line'. This collaboration resulted in three poems in the new projective plane form. We also consider what might be interesting avenues for future research, such as the incorporation of octonions in poetic form.

Keywords: projective geometry; Fano plane; poetry; poetic form; line; poetic line; octonion

*The authors report there are no competing interests to declare.*


> Only by the form, the pattern,
> Can words or [maths] reach
> The stillness…
>
> <div align="right">T.S. Eliot, Burnt Norton</div>

Once, to see if it could be done, we tried to invent a poetic form based on the mathematics of projective geometry. Poetry, as geometry, is intimately concerned with shape and form. But more than this connection, to us, projective geometry seems inherently poetic, requiring as it does an imaginative leap in questioning a maxim most people consider fundamentally unquestionable: the meeting of parallel lines. Thus, it satisfies the poetic condition 'what if'.



## 1.     On finite projective geometry

There are many approaches to finite projective geometry. One could take an algebraic approach and consider the points of a finite projective space as equivalence classes of one-dimensional subspaces of a vector space over a finite field (Hirschfeld, 1998). Or one could take an axiomatic approach and define a projective geometry as an incidence structure that satisfies some conditions (see, for example, (Dembowski, 1968)). These two approaches coincide for projective spaces of dimensions higher than 2. Here we will focus on finite projective spaces of dimension 2 - finite projective planes. We will give a very brief sketch of the mathematics of finite projective planes and refer the interested reader to (Hughes & Piper, 1973) and (Hirschfeld 1998) for a more rigorous definition. Following that we will give a highly informal description for the reader whose encounter with geometry may have ceased in high school. The reader who is more interested in the poetic side could skip to section 3.

Taking an algebraic approach, the points of a finite projective plane can be taken as equivalence classes of one-dimensional subspaces of a three-dimensional vector space over a finite field. We will call this a field plane and in general the notation PG(2,q) is used, if the finite field is one with q elements. Then the 1-dimensional vector subspaces are 0-dimensional subspaces of PG(2,q), and we call these (projective) points. The 2-dimensional vector subspaces are 1-dimensional subspaces of PG(2,q) and we call them (projective) lines. The three-dimensional vector space is the whole of the projective plane. It is the incidence properties of the subspaces that are our focus here: two points lie on a unique line, and two lines meet in a unique point.



It is not easy to visualise this, perhaps. One could imagine this: take a Euclidean plane (or Cartesian plane, or XY-plane, as they might be called), which one can picture as the plane with the x-axis and y-axis. There are lots of lines. Some of them intersect in a point. Some don't - these are parallel lines. We take a line and all the lines parallel to it (all the blue lines, say, including those that are not drawn), and we say they meet at a "point at infinity", much like rail tracks appear to meet at the horizon. We do that for all the lines (so do that for all the red lines, for example, and so on). We gather all these "points at infinity" and say they all lie on a special "line at infinity". Then we have the property that every pair of points lie on a unique line, and every pair of lines meet in a unique point.

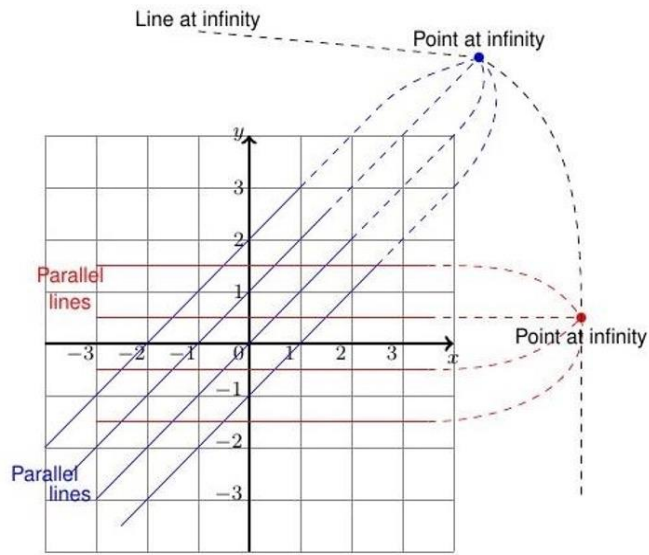

**Figure 1: An intuition of projective planes.**



(Note that when we talk about finite projective planes, the x-axis and y-axis don't go on forever, but rather wrap round, with the points and lines at infinity outside of the x- and y-axis. That can be quite hard to draw. The picture only gives some intuition to what a projective plane might look like.)

This leads us to an axiomatic approach: We might define a projective plane as a set of points and lines, where a line is a collection of points. We say a point lies on a line if that line is a collection of points that includes that point. A projective plane is such a set of points and lines that satisfies the three conditions:

1. Any two distinct points lie on a unique line.

2. Any two lines meet in a unique point.

3. There exist four points such that no three lie on one line.

It turns out that a field plane always satisfies the three conditions, but a structure that satisfies these conditions is not necessarily a field plane. Much work has been done in this area - What are those projective planes that are not field planes? What extra conditions have to be satisfied to make them field planes? See, for example, (Hughes & Piper, 1973), (Casse, 2006). On the other hand, one could ask, if we take just these three conditions, what can we say about a projective plane? There are surprisingly many structural results one could derive from them - we refer the interested reader again to (Hughes & Piper, 1973) and (Casse, 2006). For example, it transpires that a finite projective plane has the same number of lines as points, and every point lies on the same number of lines, and every line contains the same number of points. Naturally one also asks, what is the smallest set of points and lines that a projective plane can have? This gives us the Fano plane (see Figure 2): there are 7 points and 7 lines. Each point lies on 3 lines



and each line has 3 points. Every pair of lines meet in a point and every pair of points lie on a line.

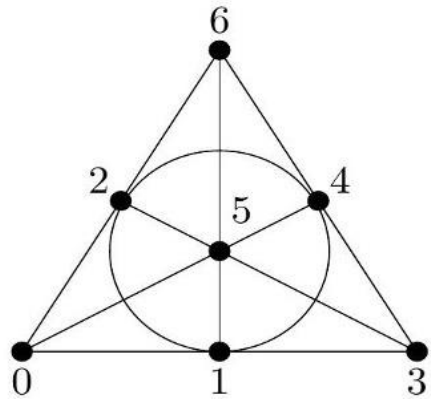

**Figure 2: The Fano plane**

The numbering of the points in this picture of the Fano plane is not accidental. If you take the numbers on one line, say, {0, 1, 3}, and take every pair of numbers and subtract one from the other modulo 7 ("modulo 7" or "mod 7" means you do your addition and subtraction like on a clock with 7 numbers 0, 1, 2,3, 4, 5, 6) you get 1,2,3,4,5,6 once each:

0 - 1 = 6 mod 7,

0 - 3 = 4 mod 7,

1 - 0 = 1 mod 7,

1 - 3 = 5 mod 7,

3 - 0= 3 mod 7,

3 - 1 = 2 mod 7.



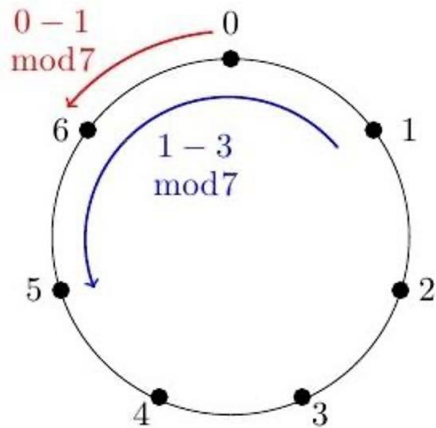

**Figure 3: Subtraction modulo 7.**

We call {0,1,3} a difference set in the set of integers modulo 7. This arises from a collineation group of PG(2,q) called the Singer group that permutes the points of PG(2,q) in a single cycle. A collineation group permutes points and lines of the plane while preserving incidence. We refer the reader to (Hughes & Piper, 1973) for details. Indeed one can define difference sets in general in a set of integers modulo n for some positive integer n as a subset that has the property that the differences modulo n between all the pairs of numbers in the subset should contain all the numbers {1, 2, … , n-1} exactly once. The existence of difference sets is related to the existence of certain types of finite projective planes.

To create our poetic form, we take the lines through the point 0, then through the point 1, and then through the point 3, to give us

0 1 3

0 4 5



0 2 6

1 5 6

1 4 2

3 4 6

3 5 2

**Figure 4: Points and lines for projective plane poetic form**

Please refer to section 3 for further details on how we translated these mathematical points and lines into the projective plane form.

**2.      On ovals and conic and algebraic forms**

The poem 'The smallest finite projective plane 4' is a found poem (we explain this later in section 3) with words found in the celebrated paper of Beniamino Segre. Here we give a flavour of that remarkable result about ovals and conics in finite field planes of odd orders. We refer the reader who is interested in the technical aspect of this result to Segre (1955), and there is also a good exposition in (Hughes & Piper, 1197). As we remarked before, a finite projective plane has the property that every line contains the same number of points. The order of the plane is this number minus one. So, in a field plane of odd order every line has an even number of points, and every point lies on an even number of lines. In a field plane we can have points that satisfy algebraic forms - functions that take points on a plane as input and give some output - we want points that give an output of 0. We call the points that satisfy an algebraic form of degree 2 (with the property of non-singularity) a conic. A conic has the same number of points as a line and has



the property that every line of the plane contains at most two points of a conic. This follows from the property of the algebraic form, similar to the familiar quadratic equations of the Euclidean plane, where a straight line crosses the graph of the function at most twice.

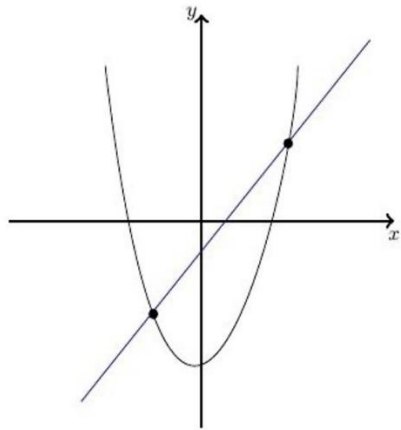

**Figure 5: Graph of a quadratic function.**

Now, without any reference to any algebra, we can also define a set of points of a finite projective plane (not necessarily the field plane) that has the property that every line of the plane contains at most two points of this set. We call this an oval. An oval is a purely combinatorial object, only concerned with incidence properties of points and lines. However, the result of (Segre, 1955) proves that in fact, in a field plane of odd order, ovals are conics and conics are ovals. It links the algebraic and the combinatorial aspects of an object in a projective plane beautifully. (In only three pages!).

The application of the first two axioms for projective planes to literature is not new. In (Hart, 2023) some examples are given. For example, the 20th century French novelist and poet Raymond Queneau proposed that one could have a literary form where words and sentences



corresponds to points and lines to satisfy the axioms. Hart herself composed a beautiful concise 'Fano fiction' out of 7 words in this literary form (Hart, 2023). But we do not know of longer poetic or literary works of this form. Here we allow projective points to correspond to lines of a poem, and projective lines to stanzas of a poem.

## 3. On making our form

At its most straightforward, poetic form can be understood as the rules that govern the structure of a poem. These rules may specify the number of lines as well as their sequence and length, which might be measured in metric feet or syllables. There can also be rules for rhyme schemes, repetitions, or the distribution of particular sounds or letters (Chivers, 2012). One well-known example of a poetic form with a long literary history is the sonnet, which has 14 lines, and conventionally reflects upon one idea or sentiment with a 'volta' or 'turn' as part of its conclusion. The original Petrarchan sonnet is divided into two sections of eight (an octave) and six lines (a sestet) and has a rhyme scheme ABBAABBA for the octave and CDCDCD or CDECDE for the sestet, though there are many variations. The Poetry Foundation has a useful explanation of the different types of sonnet, along with examples.

The poet Don Paterson, in his lengthy reflection on the poem (2018), summarises three possible theories about what form does, based on the work of linguist Richard Cureton. The first argument is that form helps reinforce a poem's content; in other words, its function is rhetorical. The second argument is, conversely, that form deliberately disrupts meaning and in so doing draws attention to 'the very failure of normal language to provide adequate symbols for our contemporary experience' (p. 345). The third argument, one which Don Paterson finds most



persuasive, as do we, is that form and content work together to create a poem's effect. Thus, while 'form' and 'content' can be distinguished from one another for the purposes of literary analysis, they 'provide the context for understanding each other, so that the poem's meaning emerges from their mutual transformation' (Hurley and O'Neill, 2012, p. 2). Ruth Padel's poem 'Revelation', for example, takes as its subject and inspiration for its form the helix structure of DNA. The poem is a 'double sonnet' of 29 lines, where line 15 marks the central point of the form, and also the point where the poem takes up the theme of the backwards running strand of DNA via an image of the caduceus (one or two snakes twirling up a staff, often a symbol of medicine or healing, based in Greek mythology). In line 15, 'But one snake, the lagging strand,/is upside down.' There are, we should note, more complex theoretical engagements with poetic form, such as Robert Shepphard's (2016, p. 1) contention that '[p]oetry is the investigation of complex contemporary realities through the means (meanings) of form'. Another interpretation is that of poet Denise Levertov, whose essay on 'organic form', published in the magazine in *Poetry* in 1965, defined form as 'a method of apperception, i.e., of recognizing what we perceive, and is based on an intuition of an order, a form beyond forms, in which forms partake' (p. 420).

As we explained earlier, the first step in creating our poetic form was the sequence of numbers based on the Fano plane, ordered with a difference set in the set of integers modulo 7 (see Figure 4). There were different ways this sequence could be interpreted poetically. For our first formal experiments, we did not want to be too distracted by content, and therefore decided to adopt a method called 'found poetry'. To make a found poem, the poet takes fragments of text they have 'found' elsewhere, which could be prose texts, other poems, historical documents, road signs, adverts, and many other sources, and arranges them into metrical lines. A classic example



of found poetry is Ezra Pound's *Cantos*. A contemporary example is the book-length poem *Zong!* by M. NourbeSe Philip, made from words found in the legal case 'Gregson v Gilbert' in 1783, which determined that the Zong massacre, in which more than 130 enslaved African people were killed, was legal. Our source was the paper 'Ovals in a finite projective plane' published by Beniamino Segre in the *Canadian Journal of Mathematics* in 1955.

Briefly, we experimented with concrete poetry, in which the shape of the poem on the page represents some aspect of its content, as with 17th Century poet George Herbert's poems 'Easter Wings' and 'The Altar', which take the shape of an angel's wings and an altar, respectively. We tried the shape of the Fano plane (see Figure 2), but it was obvious that, to work as a form, words, or even single letters, would need to represent the mathematical points, leading to an extremely sparse form. Such forms, often inspired by Japanese forms such as the haiku, tanka, or sijo, are available and can result in beautiful poems. And there are more contemporary experiments such as sonnets with one word per line, and 'fibs', which are short forms based on the Fibonacci sequence. A spare concrete form based on the Fano plane could work, but we found it creatively insufficient. The form we found most satisfying takes each mathematical point as a poetic line, and each mathematical line as a stanza (a grouping of lines, which often marks a shift in mood, time, or thought in modern free verse), which results in a form somewhat like a villanelle. The villanelle is based on a set number of lines that repeat in a particular sequence, famous examples include Dylan Thomas's 'Do Not Go Gentle into That Good Night' and Elizabeth Bishop's 'One Art'. Our repetitions seem more haphazard than the traditional villanelle, though of course they are not random. Our first formally successful poem, found in Segre's paper, was 'The smallest finite projective plane 4':



**The smallest finite projective plane 4**

It is well known that every oval can be
deemed "implausible" and
this result fills a gap in the finite.

It is well known that every oval can be
a connection between the world and
every inscribed triangle.

It is well known that every oval can be
any given oval, an arbitrary point
and therefore, a centre of perspective

deemed "implausible".
Therefore, a centre of perspective –
a connection between the world,

every inscribed triangle, and
any given oval – is an arbitrary point,
deemed "implausible".

Every inscribed triangle
fills a gap in the finite and is,
therefore, a centre of perspective.

This result fills a gap in the finite:
a connection between the world and
any given oval. An arbitrary point.



As we have previously explained poetic form as 'rules that govern the structure', it is important to note that in addition to the order of the poetic lines, there were other rules we used to make 'The smallest finite projective plane 4'. One was that within a line, the found text should remain in the same sequence as the original paper. The second was that punctuation and conjunctions (words to connect clauses such as 'and', 'the', 'as') could be adjusted in service of the poem's overall flow.  The third was that poetic lines could be ended at grammatical points or enjambed. This latter point perhaps merits further explanation. In some forms, as mentioned previously, the length of a line is predetermined. Examples include the haiku, in which line lengths are determined by syllables; another is iambic pentameter, a way of determining line length using the rhythm of poetic metre, used in many classic English poems and dramas. There are of course many more in languages other than English. In our case, line length is not determined by these factors but, as is common in contemporary poetry, by the poet's aesthetic judgement. In our form, lines may be 'end stopped' in a way that is closer to the rhythm of prose or speech, whether punctuation is used or not. An example is the poetic line 'any given oval – is an arbitrary point'. A sentence may also run over two or more poetic lines, even stanzas, using a technique called 'enjambment'. An example of enjambed lines can be found in stanzas four and five, 'Therefore, a centre of perspective – / a connection between the world, // every inscribed triangle, and / any given oval – is an arbitrary point'.

## 4. On the dream phase

I (Katherine) read 'The smallest finite projective plane 3' at a poetry critique group. The other poets were intrigued by the form, but did not connect with the mathematical language. 'Could



you use this form to write about, say, bluebells?' another poet asked. The answer is, yes, you can. To make 'The smallest finite projective plane 5' I used the same formal rules as before, with the exception that the content is not 'found' in another text, but original.

Readers might wonder, how does a poet think of what to write? The novelist Shelley Weiner once told me that there are different phases of creative writing. The first is the 'dream' phase, where critical faculties are suspended, and the writer is guided by instinct. In my case there is usually a thought or image that provokes something. With this poem it was the question, why do we assume a poet is more equipped to write about bluebells than mathematics, and if that is so, what can a poet say about projective geometry? So, I began there, creating free associative prose text, then breaking it into lines, placing the lines in their formal sequence, making adjustments.

      **The smallest finite projective plane 5**
      *'Could you use this form to write about, say, bluebells?'*

      What exactly can a poet say about
      bluebells that she can't tell us about
      parallel lines; and anyway

      what exactly can a poet say about
      the finite space below the arc of
      a stem, parallel to the others.

      What exactly can a poet say about
      the bulb that holds it all together and



also meets, somewhere in the air,

bluebells. What can she tell us about
a stem that's parallel to the others but
also meets, somewhere in the air,

bluebells, if she can't tell us about
the finite space below the arc of
the bulb that holds it all together, like

parallel lines? And anyway
the finite space below the arc
also meets, somewhere in the air,

parallel lines. And anyway
one stem might be parallel to the others, but
the bulb holds it all together.

In a previous iteration, there was a difficult part in this poem, stanzas four and five, where the possibilities for poetic lines either lost their sense or stretched too long for a reading breath. Still working in the dream phase, guided primarily by instinct, I tried to identify the cause of the problem by reading aloud, over and over. When this didn't help, I joined up the lines in all the poems I had written so far into prose-like sentences, to see if that might reveal what was wrong in poem 5. I found that sentence length varied between the poems - poem 5 had the fewest, poem 6 the most - but there didn't seem to be a discernible pattern that would show what was amiss. After some time thinking, reading and re-reading the lines, I finally spotted the problem with those two stanzas of poem 5: I had used the verb 'to say' in two different poetic lines. Originally,



mathematical lines 0 and 1 were 'What exactly can a poet say about / bluebells that can't be said about'. Changing line 1 to 'bluebells that *she can't tell us* about' resolved the issue. It might be tempting to make this a formal rule - don't use the same verb in different lines. But the thing about poetry is, I am certain it could be made to work, indeed made a feature of another poem, to do exactly that.

Speaking of sense and breathing there is, in my experience of writing in this form, a whiff of nonsense poetry, where the sounds and rhythm takes over from a literal interpretation of the words. It weaves between sense and almost sense. This whisper of nonsense raises its voice in the poem 'The smallest finite projective plane 3', with its text found in the Segre paper. We didn't consider this poem a success creatively, but it is quite satisfying rhythmically and therefore enjoyable to read aloud.

**The smallest finite projective plane 3**
Those permutations which preserve not only symmetries but reciprocities
are "self-dual," their symbol is abbreviated to an arc or branch,
and the vertices and edges may be a pair of tetrahedra, each a face-plane

of those permutations which preserve not only symmetries, but reciprocities
of the other, and every other symbol indicates which has value and
with reflection it is possible to select square faces topologically

those permutations which preserve not only symmetries but reciprocities
by opposites, as in the passage from spherical to elliptic space, which yields a map
while the rest are zero, joined by an edge combined with reflections which



are "self-dual," their symbol is abbreviated to an arc or branch

with reflection, it is possible to select square faces topologically

while the rest are zero, joined by an edge combined with reflections which

are "self-dual," their symbol is abbreviated to an arc or branch

of the other, and every other symbol indicates which has value

by opposites, as in the passage from spherical to elliptic space, which yields a map

and the vertices and edges may be a pair of tetrahedra, each a face-plane

of the other, and every other symbol indicates which has value

while the rest are zero, joined by an edge combined with reflections

and the vertices and edges may be a pair of tetrahedra, each a face-plane

with reflection, it is possible to select square faces topologically

by opposites, as in the passage from spherical to elliptic space which yields a map.

The most recent Fano plane poem we created was based on a fragment of poetic text about Octopus arms and the feeling I sometimes get in the dream phase of writing that the part of my brain I use while writing like this actually lives in my fingers.  The way I resolve problems in the dream phase is by reading aloud. Not reading silently and imagining the sounds, but actually making them, breathing out the poem to myself, stressing and unstressing the sounds and rhythm to feel out the sense. How similar, I wonder, is this process to mathematical creativity? By 'this process', I mean the weaving of sense, almost-sense and non-sense with the rhythms of typing fingers, of sound and breathing.

**The smallest finite projective plane 6**



Did you know that if an octopus loses an arm it regrows?

First it heals, then a knot forms

a small tendril grows and thickens.

Did you know that if an octopus loses an arm it regrows?

It is not clear whether this new form has its own personality

or if it is, in character and temperament, the replica of its former self.

Did you know that if an octopus loses an arm it regrows?

If this were to happen to one of my octopus fingers, would I notice the change

in the neurons in my eight fingers where the words live?

First it would heal, and then a knot would form,

in character and temperament a replica of its former self,

in the neurons of my eight fingers where the words live.

First it heals, then a knot forms.

It is not clear whether this new form has its own personality.

If this were to happen to one of my octopus fingers, would I notice the change?

A small tendril grows and thickens

and it is not clear whether this new form has its own personality

in the neurons of my eight fingers where the words live.

If a small tendril grows and thickens

in character and temperament, the replica of its former self

like one of my octopus fingers, would I notice the change?



## 5.     On critical faculties

There are, in my experience, at least two distinct sorts of critical engagement with poems: the first is the kind that takes place in what is usually called a 'workshop', where poets gather to critique one another's poems while they are in development, with the goal of helping the author to refine their work. The second is the kind of critique made after poems are published, by reviewers and literary scholars. The poems we have presented in this paper have been subjected to the first sort of critique, not yet the second. My position as author means I can offer insights about the poems from the 'inside', thus my perspective is necessarily founded in craft rather than literary criticism. On this basis, two things struck me about the poems we consider 'successful'. The first is that, to me, the poems seem inherently spatial in the way the form relates to the content. Or rather, the form encourages a spatial relationship to the content, which I can best describe as the experience of picking up a familiar object and examining it from a different angle, up close and far away, through the wrong end of a telescope, so different and increasingly unfamiliar facets are highlighted. Perhaps it can be explained as the experience of looking at a tree from the perspective of a human, then moving up the trunk to view the tree from above, as a bird, then diving below to look up at the tree from under its roots, as a worm. In this way, the form makes demands of the content, both the language itself and the images it conveys: that it be broken and reassembled, and yet still make a kind of sense. However, I can't separate this insight from the experience of crafting the poems; these spatial interpretations could be encouraged by the form, or simply by my relationship to it.

The second thing is, in the three poems we deemed successful, each seems to conclude with a kind of 'volta' (in the sonnet, this term means the 'turn' of the argument) that looks back



over the poem and offers a new perspective. This is most strongly apparent to me in poem 5 (bluebells) with the image of the bulb holding everything together, as in containing every aspect of the bluebell within its shape, as well as the more colloquial interpretation of staying calm in a confusing or stressful situation, as created by the apparently unpredictable repetitions. I note a similar 'volta' at the end of poem 4 (ovals) with 'An arbitrary point' as a gentle mockery of what I see as the vast changes in scale between which the poem moves, from a modest oval or triangle and the hubris of 'the finite'. Poem 3 does conclude with such a turn, one of the reasons we did not regard it as a creative success. So, I think we can say, with these poems at least, that the interesting spatial aspects of project geometry that inspired us creatively (as we described in the first paragraph) are reflected in the poetry. One might also ask, specifically in relation to poem 4 (ovals), to what extent the meaning of the source text is preserved. As we have explained, preserving the meaning of the original text is not a goal of found poetry (indeed, often the opposite is intended). We do not think poem 4 can be said to preserve the meaning of the Segre paper, but it may convey something of Segre's underlying insight that something which is of 'the real world' and something of the imagination, i.e. finite geometry, is bridged.

In conclusion, our initial motivation for this collaboration, to see if a poetic form based on the Fano plane could be created, our answer is yes, it can. And it can make poems that are aesthetically satisfying, to us at least. Through the process of exploring this question, though, we are left with a second: is this poetic engagement with projective planes, or something that can be read into it, interesting in any way to mathematicians?

## 6.     On poetry and mathematics



There is, of course, a rich history of connection between mathematics and different art forms (see Sriraman, 2021, for many examples). These connections are often multi-faceted, but in the case of mathematics and poetry it seems to us that they relate either there is some form of mathematical engagement with poetry or poetics, or with poems based on some aspect of mathematics, such as the Fibonacci poems previously mentioned and our projective plane forms. Within this latter category, several types of mathematical poetry have been described (for example, see Caleb Emmons's (2017) review of the *Bridges 2016 poetry anthology* edited by Sarah Glaz). This binary classification of poems based on mathematics, or the mathematical study of poetry excludes, of course, the more philosophical point made by a variety of thinkers, from Albert Einstein to Tahar Ben Jelloun, that mathematics itself is poetic. On whether the mathematical study of poetry can lead to advances in mathematical theory, Sarah Hart (2023) has pointed out that poetry, and literature more generally, has played a role in generating new mathematics in the past. Indian mathematicians from Pingala on, in studying Sanskrit poetic metre, were likely the first to use what we now call binary notation, Pascal's triangle, and the Fibonacci sequence (see also Plofker, 2009). Another example is the Renaissance architect Filippo Brunnelleschi's theory of linear perspective, which inspired the mathematician Jean-Victor Poncelet to develop the theory of projective geometry (Gamwell, 2015). However, these seem to be rare examples.

To explore whether our projective plane forms might be of interest to other mathematicians, Siaw-Lynn presented 'The smallest finite projective plane 4' at the Scottish Combinatorics Meeting held in University of St Andrews in May 2024. There were positive responses, and the poem sparked many interesting conversations about mathematics and the arts. A challenge was issued to incorporate octonions into a poem. In the most immediate



interpretation of octonions and their relationship with the Fano plane, we can view the octonions as imposing an order on the points in a line.  Octonions come from a different direction, as a sort of generalisation of imaginary numbers. The multiplication of basis elements of the octonion can be represented by the Fano plane with a particular order on the points of a line: in our labelling of the Fano plane the points on each line should go in this order (the lines can be in any order):

3 -> 1 -> 0

0 -> 2 -> 6

0 -> 5 -> 4

3 -> 4 -> 6

2 -> 5 -> 3

2 -> 1 -> 4

1 -> 5 -> 6

With this interpretation of how octonions might be incorporated into our projective plane form, we know that a poem can be made with this sequence of poetic lines, even though the sequence of the lines differs from our original form.  In another direction, we can view octonions as a mathematical object created by questioning the necessity of certain 'fundamental laws' of algebra, as projective geometry is created by questioning the necessity of Euclid's parallel postulate.  Once again, we think that this interpretation of octonions in the context of our work encourages imaginative and creative exploration, satisfying the poetic condition 'what if', and as such, could be an interesting direction for further research.



A different sort of question on which this collaboration has encouraged us to reflect is whether there is such a gulf between the ways that poets and mathematicians approach their research. As S.T. Saunders wrote in 1942, while poets and mathematicians may dream of different things ('a Paradise Lost [...] a geometry not Euclid's') we are, 'at times, singularly alike in [our] sense of detachment from so-called reality' (p. 98). This resonates with us; however we are not at all sure about Saunders's dichotomy that the mathematician is entirely devoted to her 'logic-spun domain' and the poet to her 'empire of fancy and emotion'. Nevertheless, throughout our collaboration we often asked one another, is this what it's like for you? Do you think in this way, too? How does your experience of creativity compare to mine? And often, as well, what exactly do you mean when you say…? We spent a fair amount of time and energy, for example, trying to understand one another's interpretation of 'a line'. We can sum up this relationship between mathematics and poetry, and our collaboration as a mathematician and poet, by comparing Don Paterson's view on attempts to theorise prosody. 'Many theorists have completely lost track of the fact that this material has been produced by *poets*, whose relationship to the raw material of poetry - language itself - is of such a fundamentally strange nature I doubt they would believe it' (Paterson, 2018, p. 347, emphasis in original). This recognition of strangeness, we suggest, is reflected in mathematical concepts such as octonions, in which the relationship between 'real life' numbers and what mathematicians play with is of a nature that is perhaps strange to those who are neither mathematicians nor poets.

Chivers, T. (2012) *Adventures in form: A compendium of poetic forms, rules, and constraints*. Penned in the Margins.

Dembowski, P. (1968). Finite geometries. Springer-Verlag.

Dıaz-Bánez, J. M. (2017). Mathematics and Flamenco: An Unexpected Partnership. The Mathematical Intelligencer, 39(3), 27-39.

Emmons, C. (2017): Bridges 2016 poetry anthology, Journal of Mathematics and the Arts, DOI: 10.1080/17513472.2016.1264263

Gamwell, L. (2015). Mathematics and Art: A Cultural History. Princeton University Press.

Hart, S. (2023). *Once Upon a Prime: The Wondrous Connections Between Mathematics and Literature* (2nd edition). HarperCollins.

Hirschfeld, J. W. P. (1998). Projective geometries over finite fields. (2nd edition). Clarendon Press. Oxford.

Hughes, D. R., Piper, F. C. (1973). Projective planes. Springer-Verlag.

Hurley, M.D. and O'Neill, M. (2012) *Poetic form: An introduction*. Cambridge University Press.

Levertov, D. (1965). Some notes on organic form. *Poetry*, 106(6), pp.420-425.

Paterson, D. (2018) *The poem: Lyric, sign, metre*. Faber & Faber.

Plofker, K. (2009). *Mathematics in India*. Princeton University Press.

Sanders, S. T. (1942). Mathematics and Poetry. National Mathematics Magazine, 17(3). http://www.jstor.org/stable/3028117
25

Foundation https://www.poetryfoundation.org/poems/46569/do-not-go-gentle-into-that-good-night

'One Art', by Elizabeth Bishop, in *The Complete Poems 1926-1979*. Farrar, Straus & Giroux. Also available via the Poetry Foundation https://www.poetryfoundation.org/poems/47536/one-art